*Type of Article* (Original Article)

# New Five Roots to Solve Quantic Equation in General Forms by Using Radical Expressions Along With New Theorems

**Yassine Larbaoui**[1]

[1] *Department of Electrical Engineering, Université Hassan 1er, Settat, Morocco*
*y.larbaoui@uhp.ac.ma*
*Yassine.larbaoui.uh1@gmail.com*
*Dated: (28 September 2022)*

*Abstract -* This paper presents new formulary solutions for quantic polynomial equations in general forms, where we present five solutions for any fifth degree polynomial equation with real coefficients, and thereby having the possibility to calculate the five roots of any quantic equation nearly simultaneously. The proposed roots for fifth degree polynomials in this paper are structured basing on new proposed solutions for fourth degree polynomial equations, which we developed in order to reduce the expression of any quantic polynomial to an expression of quartic polynomial.

*Keywords -* *new five solutions, new theorems, fifth degree polynomial, solving quantic equation.*

## 1. Introduction

Developing the expressions of formulary solutions to solve $n^{th}$ degree polynomial equations has been an enigmatic problem over hundreds of years to mathematicians, where many research attempts concluded to the impossibility of elaborating unified solution formulas for equations with degrees equal or higher than fifth degree by using radicals. However, finding breakthroughs to structure universal radical solutions is a living hope to solve quantic equations and above.

Solving quantic polynomial equations in general forms has been center of focus over centuries, because they presented a main problematic axe in many fields. Many adopted approaches to solve quantic equations majorly concluded to the impossibility of elaborating universal formulary expressions by using radicals as solutions, whereas others relied on treating resolvable quantic equations or reducing specific forms of fifth degree polynomials to other expressions with inferior degrees. However, developing new approaches where the scientific methodology is the principal goal, rather than aiming the discovery of specific expressions as solutions, may bring new results.

Even though solving quantic equations by using unified expressions is the center of focus in our research, developing new formulary solutions for fourth degree polynomial equations in general forms was always appealing to attention, because discovering new solution expressions for quartic equations might bring a breakthrough to express universal solutions for fifth degree polynomials by using radical expressions.

This article presents the results of a large research in mathematic science, where we present five formulary solutions for any fifth degree polynomial equation in general form with real coefficients. The developed formulary solutions for fifth degree polynomials are based on new formulary expressions as solutions for fourth degree equations, which we propose in this paper within a new theorem. We attack the problem of solving polynomial equations of $n^{th}$ degree, where n=4 and n=5, with a detailed manner basing on extendable logic, which enables us to structure new theorems by scaling from solving quartic equations upward solving quantic equations by using radicals.

Lodovico Ferrari discovered a principal solution for quartic polynomial equations in the sixteen-century (1540), but since his discovery required having a breakthrough to solve cubic equations, which was not published yet by Gerolamo Cardano, Lodovico Ferrari could not publish his elaborated solution for fourth degree polynomials immediately. In the book Ars Magna [1], which was published by Ferrari's mentor (Gerolamo Cardano), the discovered quartic solution by Lodovico Ferrari was published along with the cubic solution.

The proposed cubic solution by Cardano [2] for third degree polynomials under the form $x^3 + px + q = 0$ helps Ferrari's quartic solution to solve fourth degree polynomial equations by reducing them to the second degree. However, the elaborated solution by Ferrari does not directly help to, properly; define the four roots of any fourth degree polynomial simultaneously. The proposed method by Cardano for cubic equations was also the base to solve particular forms of $n^{th}$ degree polynomials; such as by giving radical expressions under the form $\sqrt{na + \sqrt{b}} + \sqrt{na - \sqrt{b}}$ for n=2,3,4, …, etc. [3].

There are also other proposed roots and elaborated methods to solve quartic equations such as Galois's method [4], Lagrange's method [5], algebraic geometry [6], Descartes's method [7] and Euler's solution [8].

Through the history of mathematics, the proposed solution by Cardano for cubic equations was always essential as base ground for further research about solving polynomial equations of fourth degree and above [9], [10], such as trying to prove that quantic equations do not accept quadratic expressions as potential solutions, which is further discussed in [11], [12], [13] and [14]. Nevertheless, in the results of this paper, the proposed five solutions for fifth degree polynomial equations are structured by using the proposed four roots for quartic equations in Theorem 1, where we include the expressions of quadratic roots as subparts of each solution.

There are also other recent researches dedicated to solve fourth degree polynomials where the used methods are based on reducing the expressions of equations [15-17], whereas other researches are relying on using algorithms and numerical analysis to find the roots of polynomial equations with degrees higher than four [17], [18].

In the published paper by Tschirnhaus [19], he proposed an innovative method to solve polynomial equation $P_n(x)$ of $n^{th}$ degree by transforming it into a reduced expression $Q_n(y)$ with fewer terms by extending the idea of Descartes; in which a polynomial of $n^{th}$ degree is reducible by removing its term in degree $(n-1)$. The projection of this method on quantic equations is more detailed in [20].

There are other papers treating resolvable quantic polynomials basing on radical expressions [21], [22], where the description of specific characteristics which determine whether a fifth degree polynomial may accept roots with radical expressions or not.

The advantage of this paper is proposing five formulary solutions for any form of fifth degree polynomial equations with real coefficients, and thereby having the possibility to calculate the five roots of any quantic equation nearly simultaneously. The developed solutions for fifth degree polynomial equations in this paper are based on new solutions for fourth degree polynomials in general forms, where we propose four formulary expressions as roots for any quartic equation.

Because the presented content in this paper is original, and there are many new proposed mathematical expressions and new formulas related in a scaling manner basing on extendable logic; every expressed formula will be proved mathematically and used to build the rest of content, and we will go through them by logical analysis and deduction basing on structured development.

This paper is structured as follow: section 2 where we propose new solutions for fourth degree polynomial equation in general form, section 3 where we propose new five formulary solutions for fifth degree polynomial equation in general form, and section 4 for conclusion.

## 2. New Four Solutions for Fourth Degree Polynomial Equation in General Form

This section presents new unified formulary solutions for fourth degree polynomial equation in general form (Eq.1) along with their proof. The used logic and expressions in the proof of this section is essential to prove the proposed theorems in section 3.

### 2.1. *First Proposed Theorem*

In this subsection, we propose a new theorem to resolve fourth degree polynomial equations that may be presented as shown in (Eq.1).

The expressions of proposed solutions are dependent on the value of $\left(8\frac{b^3}{a^3} - 32\frac{cb}{a^2} + 64\frac{d}{a}\right)$.

We are proposing four solutions for $\left(8\frac{b^3}{a^3} - 32\frac{cb}{a^2} + 64\frac{d}{a}\right) < 0$, four solutions for $\left(8\frac{b^3}{a^3} - 32\frac{cb}{a^2} + 64\frac{d}{a}\right) > 0$ and four solutions for $\left(8\frac{b^3}{a^3} - 32\frac{cb}{a^2} + 64\frac{d}{a}\right) = 0$.

The proposed solutions for $\left(8\frac{b^3}{a^3} - 32\frac{cb}{a^2} + 64\frac{d}{a}\right) < 0$, $\left(8\frac{b^3}{a^3} - 32\frac{cb}{a^2} + 64\frac{d}{a}\right) > 0$ and for $\left(8\frac{b^3}{a^3} - 32\frac{cb}{a^2} + 64\frac{d}{a}\right) = 0$ are expressed by using $y_{0,1}$ in (Eq.17), P in (Eq.6) and Q in (Eq.6).

The proof of this theorem is presented in an independent subsection, because it integrates long mathematical expressions.

*Theorem 1*

A fourth degree polynomial equation under the expressed form in (Eq.1), where coefficients belong to the group of numbers $\mathbb{R}$ and $a \neq 0$, has four solutions.

$$ax^4 + bx^3 + cx^2 + dx + e = 0 \ with \ a \neq 0 \quad (1)$$

If $\left(8\frac{b^3}{a^3} - 32\frac{cb}{a^2} + 64\frac{d}{a}\right) < 0$, and by using the expressions of $y_{0,1}$ in (Eq.17), P in (Eq.6) and Q in (Eq.6):

Solution 1: $S_{1,1}$ in expression (Eq.35);

Solution 2: $S_{1,2}$ in expression (Eq.36);

Solution 3: $S_{1,3}$ in expression (Eq.37);

Solution 4: $S_{1,4}$ in expression (Eq.38).

If $\left(8\frac{b^3}{a^3} - 32\frac{cb}{a^2} + 64\frac{d}{a}\right) > 0$, and by using the expressions of $y_{0,1}$ in (Eq.17), P in (Eq.6) and Q in (Eq.6):
Solution 1: $S_{2,1}$ in expression (Eq.39);

Solution 2: $S_{2,2}$ in expression (Eq.40);

Solution 3: $S_{2,3}$ in expression (Eq.41);

Solution 4: $S_{2,4}$ in expression (Eq.42).

If $\left(8\frac{b^3}{a^3} - 32\frac{cb}{a^2} + 64\frac{d}{a}\right) = 0$, and by using the expressions of $y_{0,1}$ in (Eq.28), P in (Eq.6) and Q in (Eq.6):
Solution 1: $S_{3,1}$ in expression (Eq.43);

Solution 2: $S_{3,2}$ in expression (Eq.44);

Solution 3: $S_{3,3}$ in expression (Eq.45);

Solution 4: $S_{3,4}$ in expression (Eq.46).

### 2.2. Proof of Theorem 1
By dividing the polynomial (Eq.1) on the coefficient $a$, we have the next form:
$$x^4 + \frac{b}{a}x^3 + \frac{c}{a}x^2 + \frac{d}{a}x + \frac{e}{a} = 0 \text{ with } a \neq 0 \quad (2)$$

We suppose that $x$ is expressed as shown in (Eq.3):
$$x = \frac{\left(\frac{-b}{a} + y\right)}{4} \quad (3)$$

We replace $x$ with supposed expression in (Eq.3) to reduce the form of presented polynomial in (Eq.2). Thereby, we have the presented expression in (Eq.4).
$$y^4 + y^2\left[-6\left(\frac{b}{a}\right)^2 + 16\frac{c}{a}\right] + y\left[8\left(\frac{b}{a}\right)^3 - 32\frac{cb}{a^2} + 64\frac{d}{a}\right] - 3\left(\frac{b}{a}\right)^4 + 16\frac{cb^2}{a^3} - 64\frac{db}{a^2} + 256\frac{e}{a} = 0 \quad (4)$$

To simplify the expression of polynomial equation in (Eq.4), we replace the expressions of used coefficients as shown in (Eq.5) where the values of those coefficients are defined in (Eq.6).
$$y^4 + Py^2 + Qy + R = 0 \quad (5)$$
$$P = -6\left(\frac{b}{a}\right)^2 + 16\frac{c}{a}; \quad Q = 8\left(\frac{b}{a}\right)^3 - 32\frac{cb}{a^2} + 64\frac{d}{a}; \quad R = -3\left(\frac{b}{a}\right)^4 + 16\frac{cb^2}{a^3} - 64\frac{db}{a^2} + 256\frac{e}{a} \quad (6)$$

To solve the shown polynomial equation in (Eq.5), we propose new expressions for the variable $y$; expressions (Eq.7) and (Eq.8):
$$\text{For } Q \leq 0: \quad y = \sqrt{y_0} + \sqrt{y_1} + \sqrt{y_2} \quad (7)$$
$$\text{For } Q \geq 0: \quad y = -\sqrt{y_0} - \sqrt{y_1} - \sqrt{y_2} \quad (8)$$

We propose the expressions (Eq.9) and (Eq.10) for $y_1$ and $y_2$ successively, in condition of $y_0 \neq 0$. These expressions of $y_1$ and $y_2$ are based on quadratic solutions.
$$y_1 = -\frac{\frac{P}{2} + y_0}{2} + \sqrt{\left(\frac{\frac{P}{2} + y_0}{2}\right)^2 - \frac{Q^2}{64 y_0}} \quad (9)$$

$$y_2 = -\frac{\frac{P}{2}+y_0}{2} - \sqrt{\left(\frac{\frac{P}{2}+y_0}{2}\right)^2 - \frac{Q^2}{64y_0}} \quad (10)$$

To reduce the expression of equation (Eq.5) and find a way to solve it, we propose the following expressions for the coefficients $P$ and $Q$:

$$-2[y_0^2 + y_1^2 + y_2^2] = P \quad (11)$$

$$\text{For } Q \leq 0: \quad -8\sqrt{y_0}\sqrt{y_1}\sqrt{y_2} = Q \quad (12)$$

$$\text{For } Q \geq 0: \quad 8\sqrt{y_0}\sqrt{y_1}\sqrt{y_2} = Q \quad (13)$$

In the following calculation, we replace the variable $y$ with the expression (Eq.7) where we suppose $Q < 0$, and we replace $P$ and $Q$ with their shown expressions in (Eq.11) and (Eq.12):

$$y^4 + Py^2 + Qy + R = -\left[\sqrt{y_0^4} + \sqrt{y_1^4} + \sqrt{y_2^4}\right] + 2\left[\sqrt{y_0^2}\sqrt{y_1^2} + \sqrt{y_0^2}\sqrt{y_2^2} + \sqrt{y_1^2}\sqrt{y_2^2}\right] + R$$

$$= -\left[\sqrt{y_0^2} + \sqrt{y_1^2} + \sqrt{y_2^2}\right]^2 + 4\left[\sqrt{y_0^2}\sqrt{y_1^2} + \sqrt{y_0^2}\sqrt{y_2^2} + \sqrt{y_1^2}\sqrt{y_2^2}\right] + R$$

$$= 4\left[-y_0\left(\frac{P}{2} + y_0\right) + \frac{Q^2}{64y_0}\right] + R - \frac{P^2}{4} = 0$$

$$-\left[\sqrt{y_0^2} + \sqrt{y_1^2} + \sqrt{y_2^2}\right]^2 + 4\left[\sqrt{y_0^2}\sqrt{y_1^2} + \sqrt{y_0^2}\sqrt{y_2^2} + \sqrt{y_1^2}\sqrt{y_2^2}\right] + R = 0 \Rightarrow y_0^3 + \frac{P}{2}y_0^2 + \frac{P^2-4R}{16}y_0 - \frac{Q^2}{64} = 0 \quad (14)$$

To solve the resulted expression in (Eq.14), we use Cardan's solution for third degree polynomial equations.

$$\text{For } w^3 + cw + d = 0: \quad w = \sqrt[3]{\frac{-d}{2} + \sqrt{\left(\frac{d}{2}\right)^2 + \left(\frac{c}{3}\right)^3}} + \sqrt[3]{\frac{-d}{2} - \sqrt{\left(\frac{d}{2}\right)^2 + \left(\frac{c}{3}\right)^3}} \quad (15)$$

For $y^3 + by^2 + cy + d = 0$, we use the form $y = \frac{-b+w}{3}$, and we suppose $D = 27d + 2b^3 - 9cb$ and $C = 9c - 3b^2$ to express the cubic solution as follow:

$$y = \frac{-b}{3} + \frac{1}{3}\sqrt[3]{-\frac{D}{2} + \sqrt{\left(\frac{D}{2}\right)^2 + \left(\frac{C}{3}\right)^3}} + \frac{1}{3}\sqrt[3]{-\frac{D}{2} - \sqrt{\left(\frac{D}{2}\right)^2 + \left(\frac{C}{3}\right)^3}} \quad (16)$$

By using the expression (Eq.16), the solutions of third degree polynomial equation shown in (Eq.14) are $y_{0,1}$ in (Eq.17), $y_{0,2}$ in (Eq.18) and $y_{0,3}$ in (Eq.19), where $P^i = \frac{P}{2}$, $R^i = \frac{-27Q^2 - 2P^3 + 72PR}{64}$ and $Q^i = -\frac{3P^2+36R}{16}$.

$$y_{0,1} = -\frac{P^i}{3} + \frac{1}{3}\sqrt[3]{-\frac{R^i}{2} + \sqrt{\left(\frac{R^i}{2}\right)^2 + \left(\frac{Q^i}{3}\right)^3}} + \frac{1}{3}\sqrt[3]{-\frac{R^i}{2} - \sqrt{\left(\frac{R^i}{2}\right)^2 + \left(\frac{Q^i}{3}\right)^3}} \quad (17)$$

In condition of $y_{0,1} \neq 0$, $y_{0,2}$ and $y_{0,3}$ are as follow:

$$y_{0,2} = -\frac{\frac{P}{2}+y_{0,1}}{2} + \sqrt{\left(\frac{\frac{P}{2}+y_{0,1}}{2}\right)^2 - \frac{Q^2}{64y_{0,1}}} \quad (18)$$

$$y_{0,3} = -\frac{\frac{P}{2}+y_{0,1}}{2} - \sqrt{\left(\frac{\frac{P}{2}+y_{0,1}}{2}\right)^2 - \frac{Q^2}{64y_{0,1}}} \quad (19)$$

We deduce that when $y_0$ takes the value $y_{0,1}$, the value of $y_{0,2}$ is equal to the shown value of $y_1$ in (Eq.9), and the value of $y_{0,3}$ is equal to the shown value of $y_2$ in (Eq.10).

There are three other possible expressions for $y$ which respect the proposition $-8\sqrt{y_0}\sqrt{y_1}\sqrt{y_2} = Q$ when $Q \leq 0$, and they give the same results of calculations toward having the shown third degree polynomial in (Eq.14). Thereby, they give the same values for roots $y_{0,1}$, $y_{0,2}$ and $y_{0,3}$. These three expressions are $y = -\sqrt{y_0} - \sqrt{y_1} + \sqrt{y_2}$, $y = -\sqrt{y_0} + \sqrt{y_1} - \sqrt{y_2}$ and $y = \sqrt{y_0} - \sqrt{y_1} - \sqrt{y_2}$.

By using the shown expressions in (Eq.6) and (Eq.17), the solutions of presented fourth degree polynomial equation in (Eq.5) when $\left(8\left(\frac{b}{a}\right)^3 - 32\frac{cb}{a^2} + 64\frac{d}{a}\right) < 0$ are as shown in (Eq.20), (Eq.21), (Eq.22) and (Eq.23).

$$\text{Solution 1:} \quad s_{1,1} = \sqrt{y_{0,1}} + \sqrt{-\frac{\frac{P}{2}+y_{0,1}}{2} + \sqrt{\left(\frac{\frac{P}{2}+y_{0,1}}{2}\right)^2 - \frac{Q^2}{64 y_{0,1}}}} + \sqrt{-\frac{\frac{P}{2}+y_{0,1}}{2} - \sqrt{\left(\frac{\frac{P}{2}+y_{0,1}}{2}\right)^2 - \frac{Q^2}{64 y_{0,1}}}} \quad (20)$$

$$\text{Solution 2:} \quad s_{1,2} = -\sqrt{y_{0,1}} - \sqrt{-\frac{\frac{P}{2}+y_{0,1}}{2} + \sqrt{\left(\frac{\frac{P}{2}+y_{0,1}}{2}\right)^2 - \frac{Q^2}{64 y_{0,1}}}} + \sqrt{-\frac{\frac{P}{2}+y_{0,1}}{2} - \sqrt{\left(\frac{\frac{P}{2}+y_{0,1}}{2}\right)^2 - \frac{Q^2}{64 y_{0,1}}}} \quad (21)$$

$$\text{Solution 3:} \quad s_{1,3} = -\sqrt{y_{0,1}} + \sqrt{-\frac{\frac{P}{2}+y_{0,1}}{2} + \sqrt{\left(\frac{\frac{P}{2}+y_{0,1}}{2}\right)^2 - \frac{Q^2}{64 y_{0,1}}}} - \sqrt{-\frac{\frac{P}{2}+y_{0,1}}{2} - \sqrt{\left(\frac{\frac{P}{2}+y_{0,1}}{2}\right)^2 - \frac{Q^2}{64 y_{0,1}}}} \quad (22)$$

$$\text{Solution 4:} \quad s_{1,4} = \sqrt{y_{0,1}} - \sqrt{-\frac{\frac{P}{2}+y_{0,1}}{2} + \sqrt{\left(\frac{\frac{P}{2}+y_{0,1}}{2}\right)^2 - \frac{Q^2}{64 y_{0,1}}}} - \sqrt{-\frac{\frac{P}{2}+y_{0,1}}{2} - \sqrt{\left(\frac{\frac{P}{2}+y_{0,1}}{2}\right)^2 - \frac{Q^2}{64 y_{0,1}}}} \quad (23)$$

There are three other possible expressions for $y$ which respect the proposition $8\sqrt{y_0}\sqrt{y_1}\sqrt{y_2} = Q$ when $Q \geq 0$, and they give the same third degree polynomial shown in (Eq.14) after calculations. These three expressions are $y = -\sqrt{y_0} + \sqrt{y_1} + \sqrt{y_2}$, $y = \sqrt{y_0} - \sqrt{y_1} + \sqrt{y_2}$ and $y = \sqrt{y_0} + \sqrt{y_1} - \sqrt{y_2}$.

By using the shown expressions in (Eq.6) and (Eq.17), the solutions of presented fourth degree polynomial equation in (Eq.5) when $\left(8\left(\frac{b}{a}\right)^3 - 32\frac{cb}{a^2} + 64\frac{d}{a}\right) > 0$ are as shown in (Eq.24), (Eq.25), (Eq.26) and (Eq.27).

$$\text{Solution 1:} \quad s_{2,1} = -\sqrt{y_{0,1}} - \sqrt{-\frac{\frac{P}{2}+y_{0,1}}{2} + \sqrt{\left(\frac{\frac{P}{2}+y_{0,1}}{2}\right)^2 - \frac{Q^2}{64 y_{0,1}}}} - \sqrt{-\frac{\frac{P}{2}+y_{0,1}}{2} - \sqrt{\left(\frac{\frac{P}{2}+y_{0,1}}{2}\right)^2 - \frac{Q^2}{64 y_{0,1}}}} \quad (24)$$

$$\text{Solution 2:} \quad s_{2,2} = -\sqrt{y_{0,1}} + \sqrt{-\frac{\frac{P}{2}+y_{0,1}}{2} + \sqrt{\left(\frac{\frac{P}{2}+y_{0,1}}{2}\right)^2 - \frac{Q^2}{64 y_{0,1}}}} + \sqrt{-\frac{\frac{P}{2}+y_{0,1}}{2} - \sqrt{\left(\frac{\frac{P}{2}+y_{0,1}}{2}\right)^2 - \frac{Q^2}{64 y_{0,1}}}} \quad (25)$$

$$\text{Solution 3:} \quad s_{2,3} = \sqrt{y_{0,1}} - \sqrt{-\frac{\frac{P}{2}+y_{0,1}}{2} + \sqrt{\left(\frac{\frac{P}{2}+y_{0,1}}{2}\right)^2 - \frac{Q^2}{64 y_{0,1}}}} + \sqrt{-\frac{\frac{P}{2}+y_{0,1}}{2} - \sqrt{\left(\frac{\frac{P}{2}+y_{0,1}}{2}\right)^2 - \frac{Q^2}{64 y_{0,1}}}} \quad (26)$$

$$\text{Solution 4:} \quad s_{2,4} = \sqrt{y_{0,1}} + \sqrt{-\frac{\frac{P}{2}+y_{0,1}}{2} + \sqrt{\left(\frac{\frac{P}{2}+y_{0,1}}{2}\right)^2 - \frac{Q^2}{64 y_{0,1}}}} - \sqrt{-\frac{\frac{P}{2}+y_{0,1}}{2} - \sqrt{\left(\frac{\frac{P}{2}+y_{0,1}}{2}\right)^2 - \frac{Q^2}{64 y_{0,1}}}} \quad (27)$$

Concerning $\left(8\left(\frac{b}{a}\right)^3 - 32\frac{cb}{a^2} + 64\frac{d}{a}\right) = 0$:

The expression of $y_{0,1}$ is as shown in (Eq.28) where $P^i = \frac{P}{2}$, $R^i = \frac{-2P^3+72PR}{64}$ and $Q^i = -\frac{3P^2+36R}{16}$, whereas $y_{0,2}$ and $y_{0,3}$ are as shown in (Eq.29) and (Eq.30).

$$y_{0,1} = \frac{-P^i}{3} + \sqrt[3]{-\frac{R^i}{2} + \sqrt{\left(\frac{R^i}{2}\right)^2 + \left(\frac{Q^i}{3}\right)^3}} + \sqrt[3]{-\frac{R^i}{2} - \sqrt{\left(\frac{R^i}{2}\right)^2 + \left(\frac{Q^i}{3}\right)^3}} \quad (28)$$

$$y_{0,2} = -\frac{\frac{P}{2}+y_{0,1}}{2} + \sqrt{\left(\frac{\frac{P}{2}+y_{0,1}}{2}\right)^2 - \frac{Q^2}{64y_{0,1}}} = 0 \ \ or \ \ y_{0,2} = -\left(\frac{P}{2}+y_{0,1}\right) \quad (29)$$

$$y_{0,3} = -\frac{\frac{P}{2}+y_{0,1}}{2} - \sqrt{\left(\frac{\frac{P}{2}+y_{0,1}}{2}\right)^2 - \frac{Q^2}{64y_{0,1}}} = -\left(\frac{P}{2}+y_{0,1}\right) or \ y_{0,3} = 0 \quad (30)$$

Because of having either of the expressions $(y_{0,2} = 0) \ or \ (y_{0,3} = 0)$, and having the intersection between the forms (Eq.7) and (Eq.8) for Q=0 $(y = \sqrt{y_0} + \sqrt{y_1} + \sqrt{y_2}$ for $Q \leq 0$ and $y = -\sqrt{y_0} - \sqrt{y_1} - \sqrt{y_2}$ for $Q \geq 0)$, there are four solutions for the polynomial equation shown in (Eq.5) and they are as shown in (Eq.31), (Eq.32), (Eq.33) and (Eq.34).

Solution 1: $\quad s_{3,1} = \sqrt{y_{0,1}} + \sqrt{-\left(\frac{P}{2}+y_{0,1}\right)} \quad (31)$

Solution 2: $\quad s_{3,2} = -\sqrt{y_{0,1}} - \sqrt{-\left(\frac{P}{2}+y_{0,1}\right)} \quad (32)$

Solution 3: $\quad s_{3,3} = -\sqrt{y_{0,1}} + \sqrt{-\left(\frac{P}{2}+y_{0,1}\right)} \quad (33)$

Solution 4: $\quad s_{3,4} = \sqrt{y_{0,1}} - \sqrt{-\left(\frac{P}{2}+y_{0,1}\right)} \quad (34)$

When we give the value $y_{0,1}$ in (Eq.17) to $y_0$, the values of $y_{0,2}$ in (Eq.18) and $y_{0,3}$ in (Eq.19) are equal to the shown values of $y_1$ in (Eq.9) and $y_3$ in (Eq.10) respectively. Thereby, even when we replace the value of $y_{0,1}$ in the expressions of proposed solutions by the values of $y_{0,2}$ or $y_{0,3}$, the results are only redundancies of proposed solutions, because the value of $P$ in the precedent expressions and in the proposed solutions is as follow:

$$\frac{P}{2} = -\left(\sqrt{y_0}^2 + \sqrt{y_1}^2 + \sqrt{y_2}^2\right) = -\left(\sqrt{y_{0,1}}^2 + \sqrt{y_{0,2}}^2 + \sqrt{y_{0,3}}^2\right).$$

In order to solve the polynomial equation shown in (Eq.2), we use the expression $x = \frac{-\frac{b}{a}+y}{4}$ where $y$ is the unknown variable in polynomial equation (Eq.5). By using expressions (Eq.6) and (Eq.17) for $\left(8\left(\frac{b}{a}\right)^3 - 32\frac{cb}{a^2} + 64\frac{d}{a}\right) < 0$, the solutions for equation (Eq.1) are as shown in (Eq.35), (Eq.36), (Eq.37) and (Eq.38).

Solution 1: $\quad S_{1,1} = -\frac{b}{4a} + \frac{1}{4}\sqrt{y_{0,1}} + \frac{1}{4}\sqrt{-\frac{\frac{P}{2}+y_{0,1}}{2} + \sqrt{\left(\frac{\frac{P}{2}+y_{0,1}}{2}\right)^2 - \frac{Q^2}{64y_{0,1}}}} + \frac{1}{4}\sqrt{-\frac{\frac{P}{2}+y_{0,1}}{2} - \sqrt{\left(\frac{\frac{P}{2}+y_{0,1}}{2}\right)^2 - \frac{Q^2}{64y_{0,1}}}} \quad (35)$

Solution 2: $\quad S_{1,2} = -\frac{b}{4a} - \frac{1}{4}\sqrt{y_{0,1}} - \frac{1}{4}\sqrt{-\frac{\frac{P}{2}+y_{0,1}}{2} + \sqrt{\left(\frac{\frac{P}{2}+y_{0,1}}{2}\right)^2 - \frac{Q^2}{64y_{0,1}}}} + \frac{1}{4}\sqrt{-\frac{\frac{P}{2}+y_{0,1}}{2} - \sqrt{\left(\frac{\frac{P}{2}+y_{0,1}}{2}\right)^2 - \frac{Q^2}{64y_{0,1}}}} \quad (36)$

Solution 3: $\quad S_{1,3} = -\frac{b}{4a} - \frac{1}{4}\sqrt{y_{0,1}} + \frac{1}{4}\sqrt{-\frac{\frac{P}{2}+y_{0,1}}{2} + \sqrt{\left(\frac{\frac{P}{2}+y_{0,1}}{2}\right)^2 - \frac{Q^2}{64y_{0,1}}}} - \frac{1}{4}\sqrt{-\frac{\frac{P}{2}+y_{0,1}}{2} - \sqrt{\left(\frac{\frac{P}{2}+y_{0,1}}{2}\right)^2 - \frac{Q^2}{64y_{0,1}}}} \quad (37)$

Solution 4: $\quad S_{1,4} = -\frac{b}{4a} + \frac{1}{4}\sqrt{y_{0,1}} - \frac{1}{4}\sqrt{-\frac{\frac{P}{2}+y_{0,1}}{2} + \sqrt{\left(\frac{\frac{P}{2}+y_{0,1}}{2}\right)^2 - \frac{Q^2}{64y_{0,1}}}} - \frac{1}{4}\sqrt{-\frac{\frac{P}{2}+y_{0,1}}{2} - \sqrt{\left(\frac{\frac{P}{2}+y_{0,1}}{2}\right)^2 - \frac{Q^2}{64y_{0,1}}}} \quad (38)$

By using the expression $x = \frac{-\frac{b}{a}+y}{4}$ while relying on shown expressions in (Eq.6) and (Eq.17) for $\left(8\left(\frac{b}{a}\right)^3 - \frac{32cb}{a^2} + \frac{64d}{a}\right) > 0$, the proposed solutions for equation (Eq.1) are as shown in (Eq.39), (Eq.40), (Eq.41) and (Eq.42).

Solution 1: $\quad S_{2,1} = -\frac{b}{4a} - \frac{1}{4}\sqrt{y_{0,1}} - \frac{1}{4}\sqrt{-\frac{\frac{P}{2}+y_{0,1}}{2} + \sqrt{\left(\frac{\frac{P}{2}+y_{0,1}}{2}\right)^2 - \frac{Q^2}{64y_{0,1}}}} - \frac{1}{4}\sqrt{-\frac{\frac{P}{2}+y_{0,1}}{2} - \sqrt{\left(\frac{\frac{P}{2}+y_{0,1}}{2}\right)^2 - \frac{Q^2}{64y_{0,1}}}}$ (39)

Solution 2: $\quad S_{2,2} = -\frac{b}{4a} - \frac{1}{4}\sqrt{y_{0,1}} + \frac{1}{4}\sqrt{-\frac{\frac{P}{2}+y_{0,1}}{2} + \sqrt{\left(\frac{\frac{P}{2}+y_{0,1}}{2}\right)^2 - \frac{Q^2}{64y_{0,1}}}} + \frac{1}{4}\sqrt{-\frac{\frac{P}{2}+y_{0,1}}{2} - \sqrt{\left(\frac{\frac{P}{2}+y_{0,1}}{2}\right)^2 - \frac{Q^2}{64y_{0,1}}}}$ (40)

Solution 3: $\quad S_{2,3} = -\frac{b}{4a} + \frac{1}{4}\sqrt{y_{0,1}} - \frac{1}{4}\sqrt{-\frac{\frac{P}{2}+y_{0,1}}{2} + \sqrt{\left(\frac{\frac{P}{2}+y_{0,1}}{2}\right)^2 - \frac{Q^2}{64y_{0,1}}}} + \frac{1}{4}\sqrt{-\frac{\frac{P}{2}+y_{0,1}}{2} - \sqrt{\left(\frac{\frac{P}{2}+y_{0,1}}{2}\right)^2 - \frac{Q^2}{64y_{0,1}}}}$ (41)

Solution 4: $\quad S_{2,4} = -\frac{b}{4a} + \frac{1}{4}\sqrt{y_{0,1}} + \frac{1}{4}\sqrt{-\frac{\frac{P}{2}+y_{0,1}}{2} + \sqrt{\left(\frac{\frac{P}{2}+y_{0,1}}{2}\right)^2 - \frac{Q^2}{64y_{0,1}}}} - \frac{1}{4}\sqrt{-\frac{\frac{P}{2}+y_{0,1}}{2} - \sqrt{\left(\frac{\frac{P}{2}+y_{0,1}}{2}\right)^2 - \frac{Q^2}{64y_{0,1}}}}$ (42)

By using the expression $x = \frac{-\frac{b}{a}+y}{4}$ while relying on expressions (Eq.6) and (Eq.28) for $\left(8\left(\frac{b}{a}\right)^3 - \frac{32cb}{a^2} + \frac{64d}{a}\right) = 0$, the proposed solutions for equation (Eq.1) are as shown in (Eq.43), (Eq.44), (Eq.45) and (Eq.46).

$$\text{Solution 1:} \quad S_{3,1} = -\frac{b}{4a} + \frac{1}{4}\sqrt{y_{0,1}} + \frac{1}{4}\sqrt{-\left(\frac{P}{2}+y_{0,1}\right)} \quad (43)$$

$$\text{Solution 2:} \quad S_{3,2} = -\frac{b}{4a} - \frac{1}{4}\sqrt{y_{0,1}} - \frac{1}{4}\sqrt{-\left(\frac{P}{2}+y_{0,1}\right)} \quad (44)$$

$$\text{Solution 3:} \quad S_{3,3} = -\frac{b}{4a} - \frac{1}{4}\sqrt{y_{0,1}} + \frac{1}{4}\sqrt{-\left(\frac{P}{2}+y_{0,1}\right)} \quad (45)$$

$$\text{Solution 4:} \quad S_{3,4} = -\frac{b}{4a} + \frac{1}{4}\sqrt{y_{0,1}} - \frac{1}{4}\sqrt{-\left(\frac{P}{2}+y_{0,1}\right)} \quad (46)$$

## 3. New Five Formulary Solutions for Quantic Polynomial Equation in general Form

In this section, we propose five new formulary solutions for fifth degree polynomial equation in general form shown in (Eq.47), where we rely on our proposed solutions for fourth degree polynomial equations to develop the structure of proposed solutions for quantic equation. We rely on used logic in precedent theorem (Theorem 1) by projection on quantic equations to prove the expressions of developed roots for fifth degree polynomials.

### 3.1. Second Proposed Theorem

In this subsection, we present our second proposed theorem to introduce new five formulary solutions for fifth degree polynomial equation in general form shown in (Eq.47). First, we divide the polynomial (Eq.47) on $A$ and then we replace $w$ with $\frac{-\frac{B}{A}+x}{5}$ to reduce the polynomial expression to the simplified form shown in (Eq.48) where coefficients are expressed as shown in (Eq.49), (Eq.50), (Eq.51) and (Eq.52).

$$Aw^5 + Bw^4 + Cw^3 + Dw^2 + Ew + F = 0 \text{ with } A \neq 0 \quad (47)$$

$$x^5 + cx^3 + dx^2 + ex + f = 0 \quad (48)$$

$$c = -10\frac{B^2}{A^2} + 25\frac{C}{A} \quad (49)$$

$$d = 20\frac{B^3}{A^3} - 75\frac{CB}{A^2} + 125\frac{D}{A} \quad (50)$$

$$e = -15\frac{B^4}{A^4} + 75\frac{CB^2}{A^3} - 250\frac{DB}{A^2} + 625\frac{E}{A} \quad (51)$$

$$f = 4\frac{B^5}{A^5} - 25\frac{CB^3}{A^4} + 125\frac{DB^2}{A^3} - 625\frac{EB}{A^2} + 3125\frac{F}{A} \quad (52)$$

$$z^4 + \Gamma_3 z^3 + \Gamma_2 z^2 + \Gamma_1 z + \Gamma_0 = 0 \quad (53)$$

**Theorem 2**

After reducing the form of fifth degree polynomial shown in (Eq.47) to the presented form in (Eq.48) where coefficients are expressed in (Eq.49), (Eq.50), (Eq.51) and (Eq.52); the fifth degree polynomial equation shown in (Eq.48), where coefficients belong to the group of numbers $\mathbb{R}$, can be reduced to a fourth degree polynomial equation, which may be expressed as shown in (Eq.53). The reduction from quantic polynomial to quartic polynomial is conducted by supposing $x = x_0 x_1 + x_0 x_2 + x_0 x_3 + x_1 x_2 + x_1 x_3 + x_2 x_3$, whereas supposing $z = (x_0 + x_1 + x_2 + x_3)$ is the solution for fourth degree polynomial equation in (Eq.53) by using Theorem 1 and relying on the expression $x_3 = -\frac{\Gamma_3}{4}$. The variable $\Gamma_3$ is the solution for the polynomial equation shown in (Eq.54), whereas the coefficients of this equation are shown in (Eq.55), (Eq.56) and (Eq.57). The coefficients $\Gamma_2, \Gamma_1$ and $\Gamma_0$ of quartic equation (Eq.53) are determined by using calculated values of $\Gamma_3$ and using the shown expressions in (Eq.58), (Eq.59) and (Eq.60). As a result, we have eight calculated values as potential solutions for fifth degree polynomial equation shown in (Eq.48), where many of them are only redundancies of others, because there are only five official solutions to determine.

The eight solutions to calculate for quantic equation (Eq.48) are as shown in the groups (Eq.89) and (Eq.90). The proposed five values as official solutions for fifth degree polynomial equation shown in (Eq.48) are as presented in (Eq.91), (Eq.92), (Eq.93), (Eq.94) and (Eq.95). The proposed five values as official solutions for fifth degree polynomial equation shown in (Eq.47) are as presented in (Eq.96), (Eq.97), (Eq.98), (Eq.99) and (Eq.100).

$$\lambda_2 \Gamma_3^4 + \lambda_1 \Gamma_3^2 + \lambda_0 = 0 \quad (54)$$

$$\lambda_2 = 1024 \frac{\left(e - \frac{c^2}{4}\right)}{16d}; \quad (55)$$

$$\lambda_1 = 512c - 40 \frac{(16d)^2}{e - \frac{c^2}{4}}; \quad (56)$$

$$\lambda_0 = -128 \frac{(16d)^2}{\left(e - \frac{c^2}{4}\right)^2} [f - \frac{cd}{2}]; \quad (57)$$

$$\Gamma_2 = \frac{(16d)^2}{2\Gamma_3^2 \left(e - \frac{c^2}{4}\right)} \quad (58)$$

$$\Gamma_1 = -\frac{1}{4}\Gamma_3^3 + \frac{3(16d)^2}{16\left(e - \frac{c^2}{4}\right)\Gamma_3} \quad (59)$$

$$\Gamma_0 = -\frac{1}{16}\Gamma_3^4 - \frac{(16d)^2}{32\left(e - \frac{c^2}{4}\right)} + \frac{\left(f - \frac{cd}{2}\right)(16d)^2}{2\Gamma_3^2 \left(e - \frac{c^2}{4}\right)^2} + \frac{(16d)^4}{16\Gamma_3^4 \left(e - \frac{c^2}{4}\right)^2} \quad (60)$$

**3.2. Proof of Theorem 2**

Considering the quantic equation shown in (Eq.48), we propose the expression (Eq.61) in order to reduce the form of quantic equation to a fourth degree polynomial equation. We also propose the expression (Eq.62), which presents the solution form of quartic equation by extending the used logic and presented solutions in Theorem 1.

$$x = x_0 x_1 + x_0 x_2 + x_0 x_3 + x_1 x_2 + x_1 x_3 + x_2 x_3 \quad (61)$$

$$z = x_0 + x_1 + x_2 + x_3 \quad (62)$$

We replace $x$ with its proposed value in (Eq.61), in order to end by calculations to the reduced form shown in (Eq.53).

In the shown expressions in (Eq.63), we rely on the use of $x_i$, $x_j$ and $x_k$ where $\{x_i, x_j, x_k\} \in \{x_0, x_1, x_2, x_3\}$ and $i \neq j \neq k$.

$$\alpha_1 = \sum_{i=0}^{i=3} x_i^2 \; ; \; \alpha_2 = \sum_{i \neq j} x_i^2 x_j^2 \; ; \; \alpha_3 = \sum_{i \neq j \neq k} x_i x_j x_k \; ; \; \alpha_4 = x_0 x_1 x_2 x_3 \quad (63)$$

$$x = \frac{[z^2 - \alpha_1]}{2} \quad (64)$$

$$x^2 = \alpha_2 + 2\alpha_3 z + 6\alpha_4 \quad (65)$$

$$x^3 = z^3 \alpha_3 + \frac{1}{2} z^2 [\alpha_2 + 6\alpha_4] - z\alpha_1 \alpha_3 - \frac{1}{2}[\alpha_2 + 6\alpha_4]\alpha_1 \quad (66)$$

$$x^5 = 2z^4 \alpha_3^2 + 2z^3 [\alpha_2 + 6\alpha_4]\alpha_3 + \frac{1}{2} z^2 [(\alpha_2 + 6\alpha_4)^2 - 4\alpha_3^2 \alpha_1] - 2z\alpha_3 \alpha_1 [\alpha_2 + 6\alpha_4] - \frac{1}{2}\alpha_1[\alpha_2 + 6\alpha_4]^2 \quad (67)$$

$$\gamma_4 z^4 + \gamma_3 z^3 + \gamma_2 z^2 + \gamma_1 z + \gamma_0 = 0 \quad (68)$$

We use the expressions of $\{\alpha_1, \alpha_2, \alpha_3, \alpha_4\}$ in (Eq.63), $x$ in (Eq.64), $x^2$ in (Eq.65), $x^3$ in (Eq.66) and $x^5$ in (Eq.67), to have the fourth degree polynomial shown in (Eq.68) where the values of coefficients are as follow:

$$\gamma_4 = 2\alpha_3^2$$

$$\gamma_3 = 2[\alpha_2 + 6\alpha_4]\alpha_3 + c\alpha_3$$

$$\gamma_2 = \frac{1}{2}[(\alpha_2 + 6\alpha_4)^2 - 4\alpha_3^2 \alpha_1] + \frac{1}{2}c[\alpha_2 + 6\alpha_4] + \frac{e}{2}$$

$$\gamma_1 = -2\alpha_3 \alpha_1 [\alpha_2 + 6\alpha_4] - c\alpha_1 \alpha_3 + 2d\alpha_3$$

$$\gamma_0 = -\frac{1}{2} \alpha_1 [\alpha_2 + 6\alpha_4]^2 - \frac{1}{2} c[\alpha_2 + 6\alpha_4]\alpha_1 + d[\alpha_2 + 6\alpha_4] - \frac{e\alpha_1}{2} + f$$

We divide the polynomial equation shown in (Eq.68) on $\gamma_4$ to simplify its expression. As a result, we have the shown expression in (Eq.53) where the values of coefficients are as follow:

$$\Gamma_3 = \frac{[\alpha_2 + 6\alpha_4] + \frac{c}{2}}{\alpha_3} \Rightarrow \alpha_2 = \alpha_3 \Gamma_3 - \frac{c}{2} - 6\alpha_4$$

$$\Gamma_2 = \frac{\Gamma_3^2}{4} - \alpha_1 + \frac{e - \frac{c^2}{4}}{4\alpha_3^2}$$

$$\Gamma_1 = -\Gamma_3 \alpha_1 + \frac{d}{\alpha_3}$$

$$\Gamma_0 = -\frac{1}{4}\alpha_1 \Gamma_3^2 + \frac{d\left[\alpha_3 \Gamma_3 - \frac{c}{2}\right]}{2\alpha_3^2} + \frac{f - \frac{e\alpha_1}{2} + \frac{c^2 \alpha_1}{8}}{2\alpha_3^2}$$

From precedent section (section 2), we ended with solutions expressed as $y = y_0 + y_1 + y_2$ for fourth degree polynomial equation in simple form shown in (Eq.5), whereas the solution for fourth degree polynomial equation in complete form is expressed as $z = -\frac{b}{4a} + \frac{1}{4}y_0 + \frac{1}{4}y_1 + \frac{1}{4}y_2$.

We replace $z$ with $\frac{-\Gamma_3 + y}{4}$, in order to reduce the form of quartic equation from expression (Eq.53) to expression (Eq.69), where the values of coefficients are as shown in (Eq.70).

$$y^4 + Py^2 + Qy + R = 0 \quad (69)$$

$$P = -6\Gamma_3^2 + 16\Gamma_2; \; Q = 8\Gamma_3^3 - 32\Gamma_2 \Gamma_3 + 64\Gamma_1; \; R = -3\Gamma_3^4 + 16\Gamma_2 \Gamma_3^2 - 64\Gamma_1 \Gamma_3 + 256\Gamma_0 \quad (70)$$

Concerning the fourth degree polynomial equation in (eq.53), where $z$ is as shown in (eq.62), the principal proposed expressions for the solutions are $z = -\frac{\Gamma_3}{4} + \frac{1}{4}\sqrt{y_0} + \frac{1}{4}\sqrt{y_1} + \frac{1}{4}\sqrt{y_2}$ when $Q \leq 0$ and $z = -\frac{\Gamma_3}{4} - \frac{1}{4}\sqrt{y_0} - \frac{1}{4}\sqrt{y_1} - \frac{1}{4}\sqrt{y_2}$ when $Q \geq 0$; where $x_3 = -\frac{\Gamma_3}{4}$, $x_0 = \pm\frac{1}{4}\sqrt{y_0}$, $x_1 = \pm\frac{1}{4}\sqrt{y_1}$ and $x_2 = \pm\frac{1}{4}\sqrt{y_2}$. These two principal expressions are sufficient to conduct the calculations of proof, and then generalize the results by using the other expressed forms of solutions in Theorem 1.

We replace $\Gamma_2$, $\Gamma_1$ and $\Gamma_0$ with their values in function of $\{\Gamma_3, \alpha_1, \alpha_3\}$, in order to have the expressions of $P$ in (Eq.71), $Q$ in (Eq.72) and $R$ in (Eq.73).

$$P = -2\Gamma_3^2 - 16\alpha_1 + 4\frac{\left(e - \frac{c^2}{4}\right)}{\alpha_3^2} \quad (71)$$

$$Q = -32\Gamma_3\alpha_1 - 8\Gamma_3 \frac{\left(e - \frac{c^2}{4}\right)}{\alpha_3^2} + 64\frac{d}{\alpha_3} \quad (72)$$

$$R = \Gamma_3^4 - 16\Gamma_3^2\alpha_1 + 4\,\Gamma_3^2 \frac{e - \frac{c^2}{4}}{\alpha_3^2} + 64d\frac{\Gamma_3}{\alpha_3} + 256\frac{f - \frac{e\alpha_1}{2} - \frac{cd}{2} + \frac{c^2\alpha_1}{8}}{2\alpha_3^2} \quad (73)$$

By using the shown expressions in (Eq.11), (Eq.12), (Eq.13) and (Eq.14) from the proof of first proposed theorem (Theorem 1), the values of $P$, $Q$ and $R$ will be as shown in (eq74), (Eq.75) and (Eq.76) successively.

$$P = -2[(4x_0)^2 + (4x_1)^2 + (4x_2)^2] \quad \Rightarrow \quad P = -32[\alpha_1 - \frac{\Gamma_3^2}{16}] \quad (74)$$

$$Q = -8(4x_0)(4x_1)(4x_2) \quad \Rightarrow \quad Q = -\frac{8(4x_0)(4x_1)(4x_2)(4x_3)}{4x_3} = 2048\frac{\alpha_4}{\Gamma_3} \quad (75)$$

$$R = [(4x_0)^2 + (4x_1)^2 + (4x_2)^2]^2 - 4[(4x_0)^2(4x_1)^2 + (4x_0)^2(4x_2)^2 + (4x_1)^2(4x_2)^2]$$

$$\Rightarrow R = 256\left[\alpha_1 - \frac{\Gamma_3^2}{16}\right]^2 - 1024\left[\alpha_2 - \frac{\Gamma_3^2}{16}\left(\alpha_1 - \frac{\Gamma_3^2}{16}\right)\right] \quad (76)$$

We have a group of four variables $\{\alpha_1, \alpha_2, \alpha_3, \alpha_4\}$, whereas we have a group of only three equations to solve $\{P, Q, R\}$ where all of them are dependent on the value of $\Gamma_3$. Thereby, the next step is about using the appropriate logic of analysis and calculation to find the value of $\Gamma_3$ while taking advantage of the fact that having a group of four variables enables us to solve four equations.

In order to reduce the expression of $R$ in (Eq.73), and find a way to determine the value of $\Gamma_3$, we suppose that $\left(4\frac{e - \frac{c^2}{4}}{\alpha_3^2}\Gamma_3^2 + 64\frac{d\Gamma_3}{\alpha_3}\right) = 0$ where $\frac{\Gamma_3}{\alpha_3} \neq 0$. As a result, we have the shown expression in (Eq.77).

$$\frac{\Gamma_3}{\alpha_3} = -16\frac{d}{e - \frac{c^2}{4}} \quad (77)$$

From precedent calculations of dividing the polynomial equation shown in (Eq.68) on $y_4$ to simplify its expression, along with using the expression of $\Gamma_3$, we have $\alpha_2 = (\alpha_3\Gamma_3 - \frac{c}{2} - 6\alpha_4)$.

We use the expression (Eq.77) and we replace $\alpha_2$ with $(\alpha_3\Gamma_3 - \frac{c}{2} - 6\alpha_4)$, in order to pass from expression (Eq.76) to expression (Eq.78).

$$R = -3\Gamma_3^4 + 32\Gamma_3^2\alpha_1 - 1024\left[-\frac{\Gamma_3^2}{16d}\left(e - \frac{c^2}{4}\right) - \frac{c}{2} - 6\alpha_4\right] + 256\alpha_1^2 \quad (78)$$

We have the resulted equation in (Eq.79) by using the expressions of $P$ in (Eq.71) and (Eq.74).

$$-4\Gamma_3^2 + 16\alpha_1 + 4\frac{e - \frac{c^2}{4}}{\alpha_3^2} = 0 \quad (79)$$

$$\alpha_1 = \frac{\Gamma_3^4 - \frac{(16d)^2}{e - \frac{c^2}{4}}}{4\Gamma_3^2} \quad (80)$$

We have the presented polynomial equation in (Eq.81) by using the expressions of $R$ in (Eq.73) and (Eq.78).

$$-4\Gamma_3^4 + 48\Gamma_3^2 \alpha_1 - 1024\left[-\frac{\Gamma_3^2}{16d}\left(e - \frac{c^2}{4}\right) - \frac{c}{2} - 6\alpha_4\right] - 256\frac{f - \frac{e\alpha_1}{2} - \frac{cd}{2} + \frac{c^2\alpha_1}{8}}{2\alpha_3^2} + 256\alpha_1^2 = 0 \quad (81)$$

We replace $\alpha_1$ with its value in (Eq.80) and we use the deduced value of $\alpha_4$ from (Eq.75) and (Eq.82), in order to pass from the expression of polynomial equation shown in (Eq.81) to the polynomial equation shown in (Eq.83) where the coefficients are structured only by using $c, d, e$ and $f$.

$$\alpha_4 = x_1 x_2 x_3 x_4 \Rightarrow 2048\alpha_4 = -32\Gamma_3^2 \alpha_1 - 8\Gamma_3^2 \frac{e - \frac{c^2}{4}}{\alpha_3^2} + 64\Gamma_3 \frac{d}{\alpha_3} \quad (82)$$

$$\lambda_2 (\Gamma_3^2)^3 + \lambda_1 (\Gamma_3^2)^2 + \lambda_0 (\Gamma_3^2) = 0 \quad (83)$$

The values of coefficients for polynomial equation shown in (Eq.83) are as shown in (Eq.55), (Eq.56) and (Eq.57).

Since we supposed that $(4\Gamma_3^2 \frac{e - \frac{c^2}{4}}{\alpha_3^2} + 64d \frac{\Gamma_3}{\alpha_3}) = 0$, where we adopted the shown value in (Eq.77) for $\frac{\Gamma_3}{\alpha_3}$, we eliminate the root zero as solution for polynomial equation (Eq.83) and we use the quadratic solution to solve the polynomial equation shown in (Eq.54), because all coefficients are expressed only in function of $c, d, e$ and $f$. As a result, we have four possible values for $\Gamma_3$ as solutions for polynomial equation shown in (Eq.54).

We suppose that $G_{\{\Gamma_3\}}$ is the group of solutions for shown equation in (Eq.54), where these solutions are expressed as $\Gamma_{3,i}$ and $-\Gamma_{3,i}$ with $1 \leq i \leq 2$. The group of solutions $G_{\{\Gamma_3\}}$ is determined by relying on quadratic roots.

$$G_{\{\Gamma_3\}} = \{\Gamma_{3,1}, \Gamma_{3,2}, -\Gamma_{3,1}, -\Gamma_{3,2}\} \quad (84)$$

Supposing that $P' = \frac{\lambda_1}{\lambda_2}$ and $Q' = \frac{\lambda_0}{\lambda_2}$, whereas using the expressions (Eq.55), (Eq.56) and (Eq.57); the solutions $\Gamma_{3,1}$ and $\Gamma_{3,2}$ for shown equation in (Eq.54) are as follow:

$$\Gamma_{3,1}^2 = -\frac{P'}{2} + \sqrt{\left(\frac{P'}{2}\right)^2 - 4Q'}$$

$$\Gamma_{3,2}^2 = -\frac{P'}{2} - \sqrt{\left(\frac{P'}{2}\right)^2 - 4Q'}$$

After determining the values of $\Gamma_3$ by using quadratic solutions, the following step consists of solving the polynomial equation shown in (Eq.69).

The coefficients $P$ in (Eq.71) and $R$ in (Eq.73) are dependent on $\Gamma_3^2$ and $\alpha_1$, whereas $\alpha_1$ in (Eq.80) is dependent on $\Gamma_3^2$ and the coefficient $Q$ shown in (Eq.72) is dependent on $\Gamma_3$ and $\alpha_1$. Therefore, we are going to use only $\Gamma_{3,1}$ and $\Gamma_{3,2}$ to calculate the potential values of $z$, because $\{-\Gamma_{3,1}, -\Gamma_{3,2}\}$ are going only to inverse the sign of coefficient $Q$ and thereby inversing the signs of potential values of $z$ as solutions for polynomial equation shown in (Eq.53), which will not influence the potential values of $x$ as solutions for fifth degree polynomial equation shown in (Eq.48) because $x = \frac{1}{2}(z^2 - \alpha_1)$.

We use the first proposed theorem (Theorem 1) to calculate the four solutions for fourth degree polynomial equation shown in (Eq.69) after calculating the coefficients $P$, $Q$ and $R$ for each value of $\Gamma_3$ from the group $\{\Gamma_{3,1}, \Gamma_{3,2}\}$. Thereby, we have eight values to calculate as potential solutions for the polynomial equation shown in (Eq.69).

After using first proposed theorem to solve the equation (Eq.69) for each value of $\Gamma_3$ from the group $\{\Gamma_{3,1}, \Gamma_{3,2}\}$, we have two groups of potential solutions for polynomial equation shown in (Eq.69), where each group is dependent on different value of $\Gamma_3$. We express these groups of solutions as follow: $K_{\{\Gamma_{3,1}\}}, K_{\{\Gamma_{3,2}\}}$.

$$K_{\{\Gamma_{3,1}\}} = \{S_{(\Gamma_{3,1},1)}, S_{(\Gamma_{3,1},2)}, S_{(\Gamma_{3,1},3)}, S_{(\Gamma_{3,1},4)}\} \quad (85)$$

$$K_{\{\Gamma_{3,2}\}} = \{S_{(\Gamma_{3,2},1)}, S_{(\Gamma_{3,2},2)}, S_{(\Gamma_{3,2},3)}, S_{(\Gamma_{3,2},4)}\} \quad (86)$$

Concerning the fourth degree polynomial equation shown in (Eq.53), we have two groups of solutions where each group is dependent on different value of $\Gamma_3$; as shown in (Eq.87) and (Eq.88). The values of $S_{(\Gamma_{3,i},j)}$, where $1 \leq i \leq 2$ and $1 \leq j \leq 4$, are from the expressed solutions in the groups (Eq.85) and (Eq.86).

$$M_{\{\Gamma_{3,1}\}} = \left\{-\tfrac{1}{4}\Gamma_{3,1} + \tfrac{1}{4}S_{(\Gamma_{3,1},1)}, -\tfrac{1}{4}\Gamma_{3,1} + \tfrac{1}{4}S_{(\Gamma_{3,1},2)}, -\tfrac{1}{4}\Gamma_{3,1} + \tfrac{1}{4}S_{(\Gamma_{3,1},3)}, -\tfrac{1}{4}\Gamma_{3,1} + \tfrac{1}{4}S_{(\Gamma_{3,1},4)}\right\} \quad (87)$$

$$M_{\{\Gamma_{3,2}\}} = \left\{-\tfrac{1}{4}\Gamma_{3,2} + \tfrac{1}{4}S_{(\Gamma_{3,2},1)}, -\tfrac{1}{4}\Gamma_{3,2} + \tfrac{1}{4}S_{(\Gamma_{3,2},2)}, -\tfrac{1}{4}\Gamma_{3,2} + \tfrac{1}{4}S_{(\Gamma_{3,2},3)}, -\tfrac{1}{4}\Gamma_{3,2} + \tfrac{1}{4}S_{(\Gamma_{3,2},4)}\right\} \quad (88)$$

We suppose that $S'_{(\Gamma_{3,i},j)} = (-\tfrac{1}{4}\Gamma_{3,i} + \tfrac{1}{4}S_{(\Gamma_{3,i},j)})$ where $1 \leq i \leq 2$ and $1 \leq j \leq 4$, in order to simplify the expressed values in (Eq.87) and (Eq.88). Thereby, we have two groups of values as potential solutions for fifth degree polynomial equation shown in (Eq.48). These two groups are as shown in (Eq.89) and (Eq.90) where $\alpha_{(1,\Gamma_{3,i})} = \frac{\Gamma_{3,i}^4 - (16d)^2/(e - c^2/4)}{4\Gamma_{3,i}^2}$. The expression of $\alpha_{(1,\Gamma_{3,i})}$ is an extending of the shown expression of $\alpha_1$ in (Eq.80) by changing the value of $\Gamma_3$, where $\Gamma_{3,i}$ belong to the group $\{\Gamma_{3,1}, \Gamma_{3,2}\}$.

$$N_{\{\Gamma_{3,1}\}} = \left\{\tfrac{1}{2}\left[S'^2_{(\Gamma_{3,1},1)} - \alpha_{(1,\Gamma_{3,1})}\right], \tfrac{1}{2}\left[S'^2_{(\Gamma_{3,1},2)} - \alpha_{(1,\Gamma_{3,1})}\right], \tfrac{1}{2}\left[S'^2_{(\Gamma_{3,1},3)} - \alpha_{(1,\Gamma_{3,1})}\right], \tfrac{1}{2}\left[S'^2_{(\Gamma_{3,1},4)} - \alpha_{(1,\Gamma_{3,1})}\right]\right\} \quad (89)$$

$$N_{\{\Gamma_{3,2}\}} = \left\{\tfrac{1}{2}\left[S'^2_{(\Gamma_{3,2},1)} - \alpha_{(1,\Gamma_{3,2})}\right], \tfrac{1}{2}\left[S'^2_{(\Gamma_{3,2},2)} - \alpha_{(1,\Gamma_{3,2})}\right], \tfrac{1}{2}\left[S'^2_{(\Gamma_{3,2},3)} - \alpha_{(1,\Gamma_{3,2})}\right], \tfrac{1}{2}\left[S'^2_{(\Gamma_{3,2},4)} - \alpha_{(1,\Gamma_{3,2})}\right]\right\} \quad (90)$$

We have eight values as potential solutions for the fifth degree polynomial equation shown in (Eq.48). However, many of them are only redundancies of others and there are only five official solutions to determine. The variables $\{\Gamma_{3,1}, \Gamma_{3,2}\}$ are the responsible of solution redundancies from one group to other.

In order to avoid the complications of calculating the eight values from the groups $N_{\{\Gamma_{3,1}\}}$ and $N_{\{\Gamma_{3,2}\}}$, and then determining the five solutions for quantic equation shown in (Eq.48), we propose the five expressed values in (Eq.91), (Eq.92), (Eq.93), (Eq.94) and (Eq.95) as the five official solutions for fifth degree polynomial equation shown in (Eq.48).

The first four proposed values as solutions for fifth degree polynomial equation shown in (Eq.48) are from the group $N_{\{\Gamma_{3,1}\}}$, whereas the fifth value is expressed by deduction. The useless redundancies of solutions are from one group to other; therefore, we choose the first four solutions from the same group $N_{\{\Gamma_{3,1}\}}$.

In our five proposed solutions, we use the values $S'_{(\Gamma_{3,1},j)} = -\tfrac{1}{4}\Gamma_{3,1} + \tfrac{1}{4}S_{(\Gamma_{3,1},j)}$ where $1 \leq j \leq 4$ and $S_{(\Gamma_{3,1},j)}$ from $K_{\{\Gamma_{3,1}\}}$. The variable $\alpha_{(1,\Gamma_{3,1})}$ is expressed as $\alpha_{(1,\Gamma_{3,1})} = \left[\Gamma_{3,1}^4 - \frac{(16d)^2}{e - \frac{c^2}{4}}\right]/[4\Gamma_{3,1}^2]$, whereas $\Gamma_{3,1}$ is from the group $G_{\{\Gamma_3\}}$ shown in (Eq.84), which contains the solutions for polynomial equation (Eq.54). The values of $S'_{(\Gamma_{3,1},j)}$, where $1 \leq j \leq 4$, are the solutions for quartic equation (Eq.53) and they are determined by using Theorem 1.

$$S_1 = \tfrac{1}{2}\left[S'^2_{(\Gamma_{3,1},1)} - \alpha_{(1,\Gamma_{3,1})}\right] \quad (91)$$

$$S_2 = \tfrac{1}{2}\left[S'^2_{(\Gamma_{3,1},2)} - \alpha_{(1,\Gamma_{3,1})}\right] \quad (92)$$

$$S_3 = \tfrac{1}{2}\left[S'^2_{(\Gamma_{3,1},3)} - \alpha_{(1,\Gamma_{3,1})}\right] \quad (93)$$

$$S_4 = \tfrac{1}{2}\left[S'^2_{(\Gamma_{3,1},4)} - \alpha_{(1,\Gamma_{3,1})}\right] \quad (94)$$

$$S_5 = \frac{f}{S_1 S_2 S_3 S_4} \quad (95)$$

The proposed five solutions for fifth degree polynomial equation in complete form shown in (Eq.47) are as expressed in (Eq.96), (Eq.97), (Eq.98), (Eq.99) and (Eq.100).

Solution 1: $\quad S_1 = -\tfrac{B}{5A} + \tfrac{1}{10}\left[S'^2_{(\Gamma_{3,1},1)} - \alpha_{(1,\Gamma_{3,1})}\right] \quad (96)$

Solution 2: $\quad S_2 = -\tfrac{B}{5A} + \tfrac{1}{10}\left[S'^2_{(\Gamma_{3,1},2)} - \alpha_{(1,\Gamma_{3,1})}\right] \quad (97)$

Solution 3: $\quad S_3 = -\tfrac{B}{5A} + \tfrac{1}{10}\left[S'^2_{(\Gamma_{3,1},3)} - \alpha_{(1,\Gamma_{3,1})}\right] \quad (98)$

Solution 4: $\quad S_4 = -\tfrac{B}{5A} + \tfrac{1}{10}\left[S'^2_{(\Gamma_{3,1},4)} - \alpha_{(1,\Gamma_{3,1})}\right] \quad (99)$

$$\text{Solution 5}: \quad S_5 = \frac{F}{AS_1 S_2 S_3 S_4} \quad (100)$$

### 3.3. Third Proposed Theorem

In this subsection, we present a third theorem to propose new five formulary solutions for fifth degree polynomial equation in general form shown in (Eq.101), where coefficients belong to the group of numbers $\mathbb{R}$, without simplifying the polynomial by eliminating the fourth degree part. This third proposed theorem is based on the same logic and calculations of Theorem 2. However, it is distinguished by not using the expression $x = \frac{-b+y}{5}$ which we used in Theorem 2 to eliminate the fourth degree part of quantic polynomial.

The shown equation in (Eq.101) is resulted by dividing the polynomial $Ax^5 + Bx^4 + Cx^3 + Dx^2 + Ex + F$ on $A$ whereas $A \neq 0$.

$$x^5 + bx^4 + cx^3 + dx^2 + ex + f = 0 \quad (101)$$

$$b = \frac{B}{A}; c = \frac{C}{A}; d = \frac{D}{A}; e = \frac{E}{A}; f = \frac{F}{A}.$$

**Theorem 3**

The fifth degree polynomial equation shown in (Eq.101) is reducible to the quartic equation shown in (Eq.102), where coefficients belong to the group of numbers $\mathbb{R}$ without the need to eliminate the fourth degree part. The reduction from fifth degree to fourth degree is conducted by supposing $x = x_0 x_1 + x_0 x_2 + x_0 x_3 + x_1 x_2 + x_1 x_3 + x_2 x_3$, whereas supposing $z = (x_0 + x_1 + x_2 + x_3)$ is the solution for fourth degree polynomial equation shown in (Eq.102) by using Theorem 1. The variable $Y_3$ is the solution for the polynomial equation shown in (Eq.106) by using the expression of quadratic solution, whereas $x_3 = -\frac{Y_3}{4}$. The coefficients of shown polynomial in (Eq.106) are as expressed in (Eq.107), (Eq.108) and (Eq.109). The value of $Y_3$ is equal to $Y_{3,1}$, which is presented in (Eq.120). The coefficients $Y_2$, $Y_1$ and $Y_0$ are determined by using calculated value of $Y_3$ in (Eq.120) and using the shown expressions in (Eq.103), (Eq.104) and (Eq.105). The five proposed solutions for polynomial equation (Eq.101) are as shown in (Eq.121), (Eq.122), (Eq.123), (Eq.124) and (Eq.125).

$$z^4 + Y_3 z^3 + Y_2 z^2 + Y_1 z + Y_0 = 0 \quad (102)$$

$$Y_2 = \frac{[16(d-bc)]^2}{2Y_3^2 \left(e - \frac{c^2}{4}\right)} + 4b \quad (103)$$

$$Y_1 = -\frac{1}{4} Y_3^3 + \frac{3[16(d-bc)]^2}{16\left(e - \frac{c^2}{4}\right) Y_3} + 4bY_3 \quad (104)$$

$$Y_0 = -\frac{1}{16} Y_3^4 + b Y_3^2 - \frac{[16(d-bc)]^2}{32\left(e-\frac{c^2}{4}\right)} + \frac{\left(f - \frac{cd}{2} + \frac{bc^2}{4}\right)[16(d-bc)]^2}{2Y_3^2\left(e-\frac{c^2}{4}\right)^2} + \frac{b[16(d-bc)]^2}{2Y_3^2\left(e-\frac{c^2}{4}\right)} + \frac{[16(d-bc)]^4}{16 Y_3^4 \left(e-\frac{c^2}{4}\right)^2} \quad (105)$$

$$\beta_2 Y_3^4 + \beta_1 Y_3^2 + \beta_0 = 0 \quad (106)$$

$$\beta_2 = 1024 \frac{\left(e - \frac{c^2}{4}\right)}{16(d-bc)} \quad (107)$$

$$\beta_1 = 512c - 40 \frac{[16(d-bc)]^2}{e - \frac{c^2}{4}} + 1024 b^2 \quad (108)$$

$$\beta_0 = -128 \frac{[16(d-bc)]^2}{\left(e-\frac{c^2}{4}\right)^2} \left[f - \frac{cd}{2} + \frac{bc^2}{4}\right] + 128b \frac{[16(d-bc)]^2}{\left(e-\frac{c^2}{4}\right)} \quad (109)$$

### 3.4. Proof of Theorem 3

Yassine Larbaoui, Université Hassan 1er, Settat, Morocco, 1-17, September, 2022

Considering the fifth degree polynomial equation shown in (Eq.101), we use the expression (Eq.61) as solution without eliminating the fourth degree part by the supposition $x = \frac{-b+y}{5}$. By using the expression (Eq.61), we reduce the fifth degree polynomial (Eq.101) to the quartic polynomial shown in (Eq.102).

$$x^4 = 4z^2 \alpha_3^2 + 4\alpha_3 z[\alpha_2 + 6\alpha_4] + [\alpha_2 + 6\alpha_4]^2 \quad (110)$$

$$v_4 z^4 + v_3 z^3 + v_2 z^2 + v_1 z + v_0 = 0 \quad (111)$$

We rely on the expressions of $\{\alpha_1, \alpha_2, \alpha_3, \alpha_4\}$ in (Eq.63), $x$ in (Eq.64), $x^2$ in (Eq.65), $x^3$ in (Eq.66), $x^4$ in (Eq.110) and $x^5$ in (Eq.67), to express the fourth degree polynomial shown in (Eq.111) where the values of coefficients are as follow:

$$v_4 = 2\alpha_3^2$$

$$v_3 = 2[\alpha_2 + 6\alpha_4]\alpha_3 + c\alpha_3$$

$$v_2 = \frac{1}{2}[(\alpha_2 + 6\alpha_4)^2 - 4\alpha_3^2 \alpha_1] + 4b\alpha_3^2 + \frac{1}{2}c[\alpha_2 + 6\alpha_4] + \frac{e}{2}$$

$$v_1 = -2\alpha_3 \alpha_1 [\alpha_2 + 6\alpha_4] + 4b\alpha_3[\alpha_2 + 6\alpha_4] - c\alpha_1 \alpha_3 + 2d\alpha_3$$

$$v_0 = -\frac{1}{2}\alpha_1[\alpha_2 + 6\alpha_4]^2 + b[\alpha_2 + 6\alpha_4]^2 - \frac{1}{2}c[\alpha_2 + 6\alpha_4]\alpha_1 + d[\alpha_2 + 6\alpha_4] - \frac{e\alpha_1}{2} + f$$

We have the shown fourth degree polynomial in (Eq.102) after dividing the polynomial (Eq.111) on $v_4$. The coefficients of shown polynomial in (Eq.102) are as follow:

$$Y_3 = \frac{[\alpha_2 + 6\alpha_4] + \frac{c}{2}}{\alpha_3}$$

$$Y_2 = \frac{Y_3^2}{4} - \alpha_1 + 2b + \frac{e - \frac{c^2}{4}}{4\alpha_3^2}$$

$$Y_1 = -Y_3 \alpha_1 + \frac{2b\left[\alpha_3 Y_3 - \frac{c}{2}\right]}{\alpha_3} + \frac{d}{\alpha_3}$$

$$Y_0 = -\frac{1}{4}\alpha_1 Y_3^2 + \frac{b\left[\alpha_3 Y_3 - \frac{c}{2}\right]^2}{2\alpha_3^2} + \frac{d\left[\alpha_3 Y_3 - \frac{c}{2}\right]}{2\alpha_3^2} + \frac{f - \frac{e\alpha_1}{2} + \frac{c^2 \alpha_1}{8}}{2\alpha_3^2}$$

We have the fourth degree polynomial $y^4 + My^2 + Ny + O = 0$ by replacing $z$ with $\frac{-Y_3 + y}{4}$ in the polynomial (Eq.102). The coefficients $M$, $N$ and $O$ are as expressed in (Eq.112), (Eq.113) and (Eq.114).

$$M = -2Y_3^2 - 16\alpha_1 + 32b + 4\frac{\left(e - \frac{c^2}{4}\right)}{\alpha_3^2} \quad (112)$$

$$N = -32Y_3 \alpha_1 + 64bY_3 - 8Y_3 \frac{e - \frac{c^2}{4}}{\alpha_3^2} + 64\frac{d - bc}{\alpha_3} \quad (113)$$

$$O = Y_3^4 - 16Y_3^2 \alpha_1 + 32bY_3^2 + 4Y_3^2 \frac{e - \frac{c^2}{4}}{\alpha_3^2} + 64(d - bc)\frac{Y_3}{\alpha_3} + 256\frac{f - \frac{e\alpha_1}{2} - \frac{cd}{2} + \frac{bc^2}{4} + \frac{c^2 \alpha_1}{8}}{2\alpha_3^2} \quad (114)$$

In order to reduce the expression of $O$ in (Eq.114), and find a way to determine the value of $Y_3$, we suppose that $\left(4\frac{e - \frac{c^2}{4}}{\alpha_3^2} Y_3^2 + 64(d - bc)\frac{Y_3}{\alpha_3}\right) = 0$ whereas $\frac{Y_3}{\alpha_3} \neq 0$. As a result, we have the shown expression in (Eq.115).

$$\frac{Y_3}{\alpha_3} = -16\frac{d - bc}{e - \frac{c^2}{4}} \quad (115)$$

We rely on the same used logic and processes of calculation in the proof of Theorem 2, in order to continue the proof of third proposed theorem. Thereby, the variables $\alpha_1$, $\alpha_2$ and $\alpha_4$ are as expressed in (Eq.116), (Eq.117) and (Eq.118) respectively, whereas the resulted polynomial equation to determine the value of $Y_3$ is as shown in (Eq.119).

$$\alpha_1 = \frac{Y_3^4 - 8bY_3^2 - \frac{[16(d-bc)]^2}{e - \frac{c^2}{4}}}{4Y_3^2} \quad (116)$$

$$\alpha_2 = \alpha_3 Y_3 - \frac{c}{2} - 6\alpha_4 \quad (117)$$

$$2048\alpha_4 = -32Y_3^2\alpha_1 + 64bY_3^2 - 8Y_3^2\frac{e - \frac{c^2}{4}}{\alpha_3^2} + 64Y_3\frac{(d-bc)}{\alpha_3} \quad (118)$$

$$-4Y_3^4 + [48\alpha_1 - 32b]\Gamma_3^2 - 1024\alpha_2 - 256\frac{f - \frac{e\alpha_1}{2} - \frac{cd}{2} + \frac{bc^2}{4} + \frac{c^2\alpha_1}{8}}{2\alpha_3^2} + 256\alpha_1^2 = 0 \quad (119)$$

We use the shown value of $\frac{Y_3}{\alpha_3}$ in (Eq.115) and we replace $\alpha_1$, $\alpha_2$ and $\alpha_4$ with their shown expressions in (Eq.116), (Eq.117) and (Eq.118), in order to pass from equation (Eq.119) to polynomial expression $\beta_2(Y_3^2)^3 + \beta_1(Y_3^2)^2 + \beta_0(Y_3^2) = 0$ where coefficients are as presented in (Eq.107), (Eq.108) and (Eq.109).

Relying on the proof of Theorem 2, we calculate only the expression $Y_{3,1}$ shown in (Eq.120) as a root for expressed equation in (Eq.106), and then we determine the roots of quartic equation shown in (eq102).

To determine the roots of shown equation in (Eq.102), we start by calculating the values of $\alpha_1$, $\alpha_4$ and $\alpha_2$ by replacing the variable $Y_3$ with the value of $Y_{3,1}$, then we calculate the values of $M$, $N$ and $O$, and we finish by using Theorem 1.

Supposing that $M^i = \frac{\beta_1}{\beta_2}$ and $N^i = \frac{\beta_0}{\beta_2}$, whereas using the expressions (Eq.107), (Eq.108) and (Eq.109); the solution $Y_{3,1}$ for expressed polynomial equation in (Eq.106) is as follow:

$$\Gamma_{3,1} = \sqrt{-\frac{M^i}{2} + \sqrt{\left(\frac{M^i}{2}\right)^2 - 4N^i}} \quad (120)$$

As we mentioned in the proof of Theorem 2, the use of other roots of polynomial (Eq.106) in the quartic polynomial $y^4 + My^2 + Ny + O = 0$ to calculate the values of $M$, $N$ and $O$ will generate redundancies of roots for the quantic equation shown in (Eq.101).

We determine the group $K^i_{Y_{3,1}}$, which contains the four roots for quartic equation $y^4 + My^2 + Ny + O = 0$, by using Theorem 1.

$$K^i_{\{Y_{3,1}\}} = \{S_{(Y_{3,1},1)}, S_{(Y_{3,1},2)}, S_{(Y_{3,1},3)}, S_{(Y_{3,1},4)}\}$$

We present the group of roots for quartic equation shown in (Eq.102) as $M^i_{\{Y_{3,1}\}}$, which is determined by relying on the group $K^i_{\{Y_{3,1}\}}$.

$$M^i_{\{Y_{3,1}\}} = \left\{-\frac{1}{4}Y_{3,1} + \frac{1}{4}S_{(Y_{3,1},1)}, -\frac{1}{4}Y_{3,1} + \frac{1}{4}S_{(Y_{3,1},2)}, -\frac{1}{4}Y_{3,1} + \frac{1}{4}S_{(Y_{3,1},3)}, -\frac{1}{4}Y_{3,1} + \frac{1}{4}S_{(Y_{3,1},4)}\right\}$$

We present each root for quartic equation shown in (Eq.102) as $\xi_{(Y_{3,1},i)} = -\frac{1}{4}Y_{3,1} + \frac{1}{4}S_{(Y_{3,1},i)}$, whereas $S_{(Y_{3,1},i)}$ is from the group $K^i_{\{Y_{3,1}\}}$.

The proposed five solutions for quantic equation shown in (Eq.101) are as expressed in (Eq.121), (Eq.122), (Eq.123), (Eq.124) and (Eq.125). The expressions $\xi_{(Y_{3,1},1)}$, $\xi_{(Y_{3,1},2)}$, $\xi_{(Y_{3,1},3)}$ and $\xi_{(Y_{3,1},4)}$ present the calculated roots for quartic equation (Eq.102) by using Theorem 1. $Y_{3,1}$ is calculated by using the shown expression in (Eq.120). We use the expression (Eq.116) to calculate the value of $\alpha_{(1,Y_{3,1})}$; which has the following value:

$$\alpha_{(1,Y_{3,1})} = \frac{Y_{3,1}^4 - 8bY_{3,1}^2 - \frac{[16(d-bc)]^2}{e - \frac{c^2}{4}}}{4Y_{3,1}^2}$$

$$\text{Solution 1:} \quad S_1 = \frac{1}{2}\left[\xi^2_{(Y_{3,1},1)} - \alpha_{(1,Y_{3,1})}\right] \quad (121)$$

$$\text{Solution 2:} \quad S_2 = \frac{1}{2}\left[\xi^2_{(Y_{3,1},2)} - \alpha_{(1,Y_{3,1})}\right] \quad (122)$$

$$\text{Solution 3:} \quad S_3 = \frac{1}{2}\left[\xi^2_{(Y_{3,1},3)} - \alpha_{(1,Y_{3,1})}\right] \quad (123)$$

$$\text{Solution 4:} \quad S_4 = \frac{1}{2}\left[\xi^2_{(Y_{3,1},4)} - \alpha_{(1,Y_{3,1})}\right] \quad (124)$$

$$\text{Solution 5:} \quad S_5 = \frac{f}{S_1 S_2 S_3 S_4} \quad (125)$$

## 4. Conclusion

In the first presented theorem in this paper, we propose new four formulary solutions for any quartic polynomial equation in general form, which enable us to develop the formulary structures of five roots for any fifth degree polynomial equation in general form.

The proposed expressions as solutions in first theorem enable to calculate the four roots of any quartic polynomial equation nearly simultaneously, whereas the proposed solutions in second theorem enable to calculate the five roots of any fifth degree polynomial equation nearly simultaneously. The third proposed theorem is based on the same logic and calculations of Theorem 2 whereas proposing five roots for quantic equation; however, it is distinguished by not eliminating the fourth degree part of fifth degree polynomial by avoiding the use of the expression $x = \frac{-b+y}{5}$.

The second and third proposed theorems in this paper are based on reducing any fifth degree polynomial to a fourth degree polynomial, and then using Theorem 1 to solve the resulted quartic equation.

The principal criteria of presented theorems in this paper is proposing new radical solutions for quartic equations and quantic equations to enable the calculation of all roots of these equations nearly simultaneously, whereas using the expressions of quadratic roots and cubic roots as subparts of each proposed solution.